\documentclass{siamonline171218}
\usepackage[T1]{fontenc}
\usepackage[latin9]{inputenc}
\usepackage{geometry}
\usepackage{xcolor}
\usepackage{pdfcolmk}
\usepackage{verbatim}

\usepackage{lipsum}
\usepackage{amsfonts}
\usepackage{graphicx}
\usepackage{epstopdf}
\usepackage{algorithmic}
\ifpdf
  \DeclareGraphicsExtensions{.eps,.pdf,.png,.jpg}
\else
  \DeclareGraphicsExtensions{.eps}
\fi

\newsiamremark{remark}{Remark}
\newsiamremark{hypothesis}{Hypothesis}
\crefname{hypothesis}{Hypothesis}{Hypotheses}
\newsiamthm{claim}{Claim}

\headers{Sensitivity and Bifurcation Analysis of a DAE Model for a Microbial Electrolysis Cell}{H. J. Dudley, L. Lu, Z. J. Ren, and D. M. Bortz}

\title{Sensitivity and Bifurcation Analysis of a Differential-Algebraic
Equation Model for a Microbial Electrolysis Cell}

\author{Harry J. Dudley\thanks{Department of Applied Mathematics, University of Colorado, Boulder,
CO 80309-0526}
\and Lu Lu\thanks{Department of Civil, Environmental, and Architectural Engineering,
University of Colorado, Boulder, CO 80309-0428}
\and Zhiyong Jason Ren\footnotemark[2]
\and David M. Bortz\footnotemark[1] \thanks{Corresponding author (\email{dmbortz@colorado.edu})}.}

\usepackage{amsopn}

\usepackage{pgf,tikz} 
\usepackage{pgfplots}
\usetikzlibrary{spy}
\pgfplotsset{compat=1.14}

\newtheorem{condition}{Condition}[section]

\definecolor{copper}{cmyk}{0,0.9,0.9,0.2} 
\colorlet{lightgray}{black!25} 
\colorlet{darkgray}{black!75} 

\begin{document}

\maketitle

\begin{abstract}
Microbial electrolysis cells (MECs) are a promising new technology
for producing hydrogen cheaply, efficiently, and sustainably. However, to scale up this
technology, we need a better understanding of the processes in the
devices. In this effort, we present a differential-algebraic equation
(DAE) model of a microbial electrolysis cell with an algebraic constraint
on current. We then perform sensitivity and bifurcation analysis for
the DAE system. The model can be applied either to batch-cycle MECs
or to continuous-flow MECs. We conduct differential-algebraic sensitivity analysis
after fitting simulations to current density data for a batch-cycle
MEC. The sensitivity analysis suggests which parameters have the greatest
influence on the current density at particular times during the experiment.
In particular, growth and consumption parameters for exoelectrogenic
bacteria have a strong effect prior to the peak current
density. An alternative strategy to maximizing peak current density
is maintaining a long term stable equilibrium
with non-zero current density in a continuous-flow MEC.
We characterize the minimum dilution rate required for a stable nonzero
current equilibrium and demonstrate transcritical bifurcations in
the dilution rate parameter that exchange stability between several
curves of equilibria. Specifically, increasing the dilution rate transitions
the system through three regimes where the stable equilibrium exhibits 
(i) competitive exclusion by methanogens, (ii) coexistence, and 
(iii) competitive exclusion by exolectrogens. Positive
long term current production is only feasible in the final two regimes.
These results suggest how to modify system parameters to increase
peak current density in a batch-cycle MEC or to increase the long
term current density equilibrium value in a continuous-flow MEC.

\begin{keyword}
Microbial Electrolysis Cell, Differential-Algebraic Equation, Sensitivity Analysis, Bifurcation{\small \par}
\end{keyword}
\end{abstract}

\section{Introduction}

Microbial electrolysis cells (MECs) are devices that produce hydrogen
from renewable organic matter, such as wastewater. These devices require
less energy input than water electrolysis and have greater efficiency
than fermentative hydrogen production \cite{lu_microbial_2016,liu_electrochemically_2005,rozendal_principle_2006,pinto2011}.
The technology is promising, but performance remains low for MECs
and other microbial bioelectrochemical cells, despite significant
work on the experimental design and scaling up of technology \cite{oliveira_overview_2013,ortiz-martinez_developments_2015,recio-garrido_modeling_2016}.
In Section \ref{subsec:Biological-and-Electrochemical}, we describe
the biological and electrochemical processes occurring in MECs and,
in Section \ref{subsec:Mathematical-Modeling-Background}, we discuss
mathematical models that have been used to explain MEC operation.

\subsection{Biological and Electrochemical Background\label{subsec:Biological-and-Electrochemical}}

\begin{figure}
\centering{}\resizebox {\textwidth} {!} {
\begin{tikzpicture}[spy using outlines={rectangle,lens={scale=2}, size=8cm, connect spies}]
\draw (-1,0) to [controls=+(90:0.5) and +(90:0.5)] (6.5,0);     
\draw[fill=blue!60, fill opacity=0.5] (-1,-0.5) to [controls=+(90:0.5) and +(90:0.5)] (6.5,-0.5);

\draw[fill=copper, fill opacity=0.75] (-1,2) rectangle (0,-4);     
\draw (-1.4,2.3) node {\Large Anode};          

\draw[fill=gray, fill opacity=0.75] (5.5,2) rectangle (6.5,-4);     
\draw (7.1,2.3) node {\Large Cathode};          

\draw[join = round, thick] (-0.5,2) -- (-0.5,2.5) -- (6,2.5) -- (6,2);     
\draw (2.75,2.5) node [rectangle, draw, fill=red!30] {\Large PS};     
\draw (1,3) node {$e^- \rightarrow$};     
\draw (4.5,3) node {$e^- \rightarrow$};          

\foreach \y  in {-2,-2.5,-3.5}%
\shadedraw [ball color= green] (0.2,\y) circle (0.25cm);     
\foreach \y  in {-3,-4}%
\shadedraw [ball color= blue] (0.2,\y) circle (0.25cm);     
\foreach \y  in {-3,-3.5,-4}%
\shadedraw [ball color= blue] (0.4,\y,-0.5) circle (0.25cm);

\draw (-1,0) .. controls +(-90:0.5) and +(-90:0.5) .. (6.5,0);     
\draw (-1,0) .. controls +(-90:0.5) and +(-90:0.5) .. (6.5,0)         
	-- (6.5,-0.5) .. controls +(-90:0.5) and +(-90:0.5) .. (-1,-0.5) -- (-1,0);
    
\draw[fill=blue!60, fill opacity=0.5] (-1,-0.5) .. controls +(-90:0.5) and +(-90:0.5) .. (6.5,-0.5);     
\draw[fill=blue!60, fill opacity=0.5] (-1,-0.5) .. controls +(-90:0.5) and +(-90:0.5) .. (6.5,-0.5)         
	-- (6.5,-4) .. controls +(-90:0.5) and +(-90:0.5) .. (-1,-4) -- (-1,-0.5);              

\draw (0.9,-1) node {$M_\text{ox}$};     
\draw[thick] (0.6,-1.4)+(45:.25) [yscale=2,xscale=2,<-] arc(110:250:.25);     
\draw (1,-2.4) node {$M_\text{red}$};     
\draw[thick,dashed] (1,-1.95)+(-70:.25) [yscale=2,xscale=2,<-] arc(-70:70:.25);           

\spy on (0.6,-2.5) in node at (-8,1);              

\draw[->, thick] (.6,-1.5) -- (2.25,-1.5);     
\draw (2,-1.3) node {H$^+$};           

\draw[->, thick] (.375,-1.5) -- (-.375,-1.5);     
\draw (0.25,-1.3) node {e$^-$};          

\draw (1.9,-2.8) node {Organic};     
\draw[thick] (1.2,-2.8) [yscale=3,->] arc(90:320:.25);     
\draw (1.6,-3.8) node {CO$_2$};     
\draw (2,-4.15) node {+CH$_4$};          

\draw[thick] (5.25,-2.5)+(-135:.25) [yscale=3,<-] arc(-135:135:.25);     
\draw (4.7,-2.6) node {H$_2$};      

\end{tikzpicture} 
}\caption{Illustration of a single-chamber MEC. The device consists of an anode
and a cathode. Organic substrate is fed into the device. The substrate
is oxidized by exoelectrogenic microorganisms (green spheres) in a
biofilm on the anode, producing CO$_{2}$. In the process, an intracellular
mediator, $M$, is oxidized and electrons are transferred to the anode.
Methanogenic microorganisms (blue spheres) compete for the substrate,
producing CH$_{4}$ as well as CO$_{2}$ and decreasing MEC efficiency.
At the cathode, protons from the anode combine with electrons to produce
hydrogen via a reduction reaction. Microbial electrolysis is endothermic
(positive Gibbs free energy), so an external voltage must be applied
by a power source (PS).\label{fig:Cartoon_MEC}}
\end{figure}
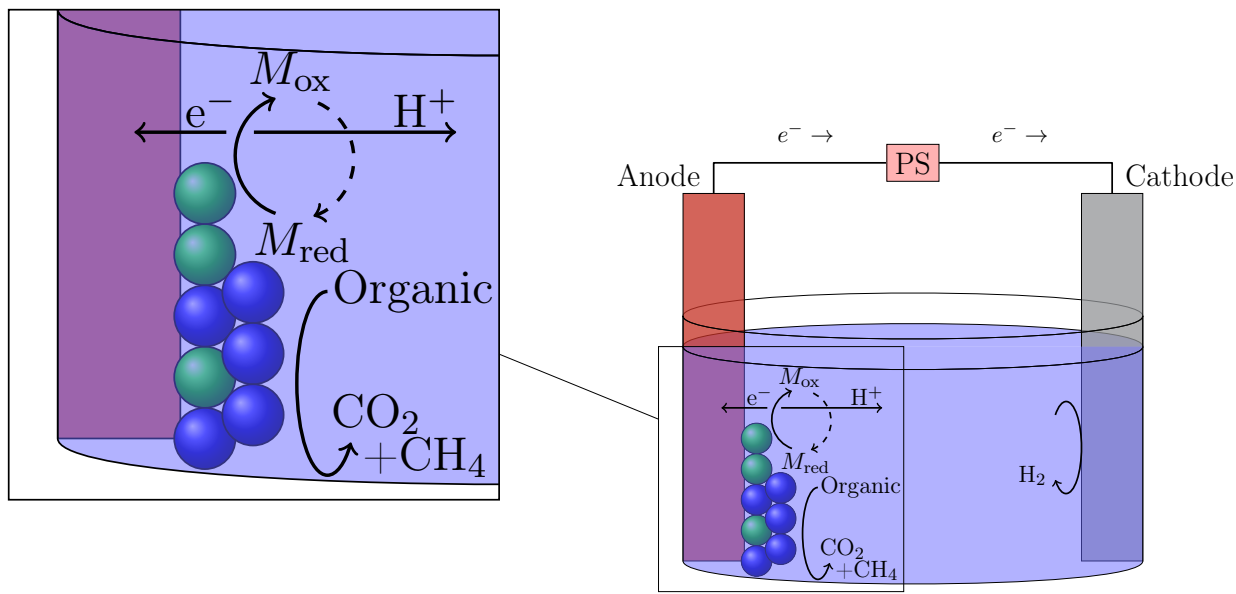

MECs are based on microbial fuel cells (MFCs). These devices employ
a biofilm of bacteria on the cell anode to biocatalyze an oxidation-reduction
reaction. The main bacteria involved are known as exoelectrogenic
microorganisms (or exoelectrogens) because they transfer electrons
extracellularly. Figure \ref{fig:Cartoon_MEC} provides a visual representation
the MEC device. The biofilm holds the microorganisms in place while
exoelectrogens (depicted by green spheres) oxidize a substrate. In
the process, electrons are transferred to the anode via either intracellular
mediators, nano-pili (or nanowires), cytochromes, or a combination
of these, depending on the specific microorganism \cite{logan_exoelectrogenic_2009}.
Current is then generated through an external circuit due to the potential
difference between anode and cathode. On the cathode side of an MEC,
the current can be used to drive a reduction reaction such as hydrogen
production \cite{liu_electrochemically_2005,rozendal_principle_2006,rabaey_microbial_2010}.
In practice, the process of microbial electrolysis is endothermic
(positive Gibbs free energy), so an external voltage must be applied.
However, the action of the exoelectrogens decreases the amount of
energy that is needed for the reaction. For the experiment described
in this paper, only 0.6 - 1.0 Volts were applied to produce hydrogen,
compared to about 1.8 - 2.0 Volts for hydrogen production via water
electrolysis \cite{lu_hydrogen_2009,lu_nickel_2016}. 

When the substrate is complex wastewater or a mixture of compounds,
fermenting microorganisms convert the complex organic matter into
simpler compounds which the exoelectrogens can consume \cite{batstone_iwa_2002,wang_comprehensive_2013}.
MEC efficiency can be decreased by a variety of factors. For example,
other types of microorganisms are sometimes introduced unintentionally
and are often present in the organic material that is fed into the
MEC. In particular, methanogenic or methane-producing microorganisms
compete with exoelectrogens for substrate, decreasing exoelectrogen
growth. (Methanogens are depicted as blue spheres in Figure \ref{fig:Cartoon_MEC}).
Methanogen activity can be measured by how much methane is produced.
A specific variety know as hydrogenotrophic methanogens can consume
some of the hydrogen produced at the cathode. Complicating matters
further, exoelectrogens themselves can consume hydrogen to accelerate
current generation while, at the same time, increasing the energy
loss on the electrodes \cite{lu_active_2016}. This consumption
of hydrogen is sometimes remedied by a two chamber design that separates
anode from cathode. 

In addition, several processes contribute to overpotentials of the
electrodes. Overpotentials are the difference in potential (or voltage)
between the observed potential and the calculated thermodynamic reduction
potential of a half reaction. In other words, overpotentials are voltage
losses or inefficiencies in the MEC. As such, they must be accounted
for in our model's current equation. MEC overpotentials include ohmic
losses, activation losses, concentration losses, and microbial metabolic
losses \cite{logan_microbial_2006}. Ohmic losses are related to
various types of resistance in the circuit. Activation losses are
related to the activation energy of the oxidation-reduction reactions
occurring in the cell. Concentration losses are caused by various
processes that limit the concentration of reactants at the anode and
the cathode. Microbial metabolic losses refer to the energy lost to
the microorganisms' metabolic pathways. A modified version of Ohm's
law that includes these voltage losses provides the algebraic constraint
in the DAE system of Section \ref{sec:Model}.

\subsection{Mathematical Modeling Background\label{subsec:Mathematical-Modeling-Background}}

Several models have been proposed to describe MFC or MEC operation,
mainly focusing on the anode reaction. Several other models describe
specific processes in MFCs. For an overview of these models, see the
review articles by Ortiz-Mart\'{i}nez et al. \cite{ortiz-martinez_developments_2015} and Recio-Garrido et al. \cite{recio-garrido_modeling_2016}. 

In 1995, Zhang and Halme presented a differential-algebraic equation
model for a MFC that used an added chemical mediator, 2-hydroxy-1,4
naphthoquinone or HNQ \cite{zhang_modelling_1995}. They used ordinary
differential equations to describe concentrations of substrate, a
reaction intermediate, and the added chemical mediator, HNQ. The model
simplified the MFC system by assuming that one type of microorganism
was present and that total biomass of the microbes was constant. However,
it did introduce several fundamental aspects of MFC modeling, including
Monod kinetics to describe substrate consumption, Faraday's law to
describe current, the Nernst equation to describe electromotive force,
and Ohm's law for ohmic overpotentials. Other overpotentials were
neglected or approximated. By 2003, researchers had discovered that
MFCs do not require an external chemical mediator. Instead, certain
microorganisms can directly transfer electrons to the anode \cite{chaudhuri_electricity_2003}.
Subsequent models may include internal mediators, but omit external
chemicals.

There are several PDE models that include biofilm growth in up to
three dimensions. Marcus et al. presented
a model of a 1D biofilm, which was modeled as a conductive solid matrix \cite{kato_marcus_conduction-based_2007}.
They also derived a Nernst-Monod relation to describe oxidation of
the electron-donor. Picioreanu et al. presented
a model that incorporated the Butler-Volmer equation to calculate
current density, allowed for multiple microorganism species including
methanogens, and explicitly modeled growth of the biofilm in one,
two, or three dimensions \cite{picioreanu_computational_2007}. 
One disadvantage of this model was that
simulating MFC operation for 15 days took about 6 minutes for a 1D
biofilm and about 14 hours for a 3D biofilm. Picioreanu et. al
later updated the model to use the International Water Association's
anaerobic digestion model, ADM1 \cite{batstone_iwa_2002}, with
six microorganism populations \cite{picioreanu_mathematical_2008}. 
The group also extended the model to
investigate effects of pH and electrode geometry \cite{picioreanu_model_2010}. 

In contrast, Pinto et al. have proposed a DAE compartment model that
simulates current in MFCs or current and hydrogen production in MECs
\cite{pinto_two-population_2010,pinto2011}. The constraint comes
from a modified version of Ohm's law that includes voltage losses
(or overpotentials) in the system. The constraint equation is transcendental
because of an inverse hyperbolic sine approximation that allows one
to use the Butler-Volmer equation for activation voltage losses. The
DAE model has the advantage of being computationally inexpensive compared
to PDE models, making it a better candidate for process control. The
major disadvantage is that biofilm is modeled as a compartment, so
biofilm properties besides population growth cannot be simulated.
However, the model does provide a reasonable description of current
density in MFCs and MECs fed on simple substrates. Pinto et al. allow
for the influent to contain a mixed organic substrate that can be
broken down into a simpler substrate such as acetate. They also consider
the action of four types of microorganisms: fermenters, methanogens,
exoelectrogens, and hydrogenotrophic methanogens. Like previous models,
Pinto et al. use multiplicative Monod kinetics for microorganism
growth. Concentration overpotential is modeled with the Nernst equation,
and activation overpotential
is modeled with an approximation to the Butler-Volmer equation, as
in previous models \cite{kato_marcus_conduction-based_2007, picioreanu_computational_2007}.

The ODE model in \cite{pinto_two-population_2010} builds upon the
classic chemostat model. The chemostat was first proposed by Monod \cite{monod_recherches_1942}
and then independently by Novick and Szilard \cite{novick_description_1950}. A simple
chemostat is a chemical reactor containing bacteria, with concentration
$X$, and a substrate for bacterial growth, with concentration $S$.
Suppose nutrients enter the reactor with concentration $S_{0}$ at
a constant, possibly zero, inflow rate, $F_{\text{in}}$. Suppose
also that fluid flows out of the reactor at the same rate. If the
reactor has volume $V$, the dilution rate is given by $D=F_{\text{in}}/V$.
We can model the reactor with following ODE system
\begin{align}
\frac{dS}{dt} & =D(S_{0}-S)-\frac{1}{Y}\frac{\mu SX}{K+S},\label{eq:chemostatSubstrate-1}\\
\frac{dX}{dt} & =\frac{\mu SX}{K+S}-DX,\label{eq:chemostatBacteria-1}
\end{align}
where $\mu$ is the maximum growth rate and $Y$ is the yield. $K$
is the half saturation of half rate constant, reflecting that growth
occurs at half the maximal rate when $K=S$. The model allows us to
quantify several features of the chemostat. For example, there is
a break-even nutrient concentration, $S=\lambda=\frac{KD}{\mu-D}$,
that is required in order for the growth rate to exceed the dilution
rate. There is also a washout equilibrium at $(S,X)=(S_{0},0)$ and
a survival equilibrium at $(S,X)=(\lambda,Y(S_{0}-\lambda))$. When
the system is perturbed slightly from the washout equilibrium, solutions
obey a linearized system with Jacobian
\[
J=\begin{bmatrix}-D & -\frac{1}{Y}\frac{\mu S_{0}}{K+S_{0}}\\
0 & \frac{\mu S_{0}}{K+S_{0}}-D
\end{bmatrix}.
\]
This tells us that the washout equilibrium is stable if $\frac{\mu S_{0}}{K+S_{0}}<D$
and unstable if $\frac{\mu S_{0}}{K+S_{0}}>D$ . Similarly, a positive
survival equilibrium exists if and only if $\frac{\mu S_{0}}{K+S_{0}}>D$. 

Additional conclusions can be reached if we consider a chemostat in
which two species of bacteria compete for the substrate. This scenario
may be modeled by the ODE system
\begin{align}
\frac{dS}{dt} & =D(S_{0}-S)-\frac{1}{Y_{1}}\frac{\mu_{1}SX_{1}}{K_{1}+S}-\frac{1}{Y_{2}}\frac{\mu_{2}SX_{2}}{K_{2}+S},\label{eq:chemostatSubstrate-2}\\
\frac{dX_{1}}{dt} & =\frac{\mu_{1}SX_{1}}{K_{1}+S}-DX_{1},\label{eq:chemostatBacteria-2a}\\
\frac{dX_{2}}{dt} & =\frac{\mu_{2}SX_{2}}{K_{2}+S}-DX_{2}.\label{eq:chemostatBacteria-2b}
\end{align}
In this scenario, the break even concentrations are $\lambda_{1}=\frac{K_{1}D}{\mu_{1}-D}$
and $\lambda_{2}=\frac{K_{2}D}{\mu_{2}-D}$. Coexistence can only
occur when $\lambda_{1}=\lambda_{2}$. In general, $\lambda_{1}<\lambda_{2}$,
meaning that species $X_{1}$ can survive at lower substrate concentrations,
although this ordering may be determined by the dilution rate, $D$.
When $\lambda_{1}<\lambda_{2}$, species $X_{1}$ will survive and
the model will tend to a stable equilibrium at $(S,X_{1},X_{2})=(\lambda_{1},Y_{1}(S_{0}-\lambda_{1}),0).$
This phenomenon is known as the principle of competitive exclusion.
Hsu provided a mathematical proof  \cite{hsu_mathematical_1977} and the
phenomenon was later tested experimentally by Hansen and Hubbell \cite{hansen_single-nutrient_1980}.
For more analysis, see the text on chemostat theory by Smith and Waltman \cite{smith_theory_1995}.

\subsection{Overview}

The rest of the document is organized as follows. Section \ref{sec:Model}
describes the experiment that produced our data as well as the MEC
model and the parameter estimation process. Differential-algebraic
sensitivity analysis is presented in Section \ref{sec:Sensitivity-Equations}.
This method provides a way of looking at sensitivity of current density
to parameters. Bifurcation analysis of the model with respect to the
dilution rate parameter is shown in Section \ref{sec:Bifurcation-Analysis}.
This provides critical information about how changing the dilution
rate can change the stability of equilibria. Finally, the model results
are discussed in Section \ref{sec:Discussion} and concluding remarks
appear in Section \ref{sec:Conclusion}.

\newpage

\section{Model and Data\label{sec:Model}}

\subsection{Experiment}

Our data were collected from single-chamber membraneless MECs with
liquid volume $V=90$ mL and carbon brush anodes with surface area
$A_{\text{sur,A}}=0.4$ meters squared. The reactors were operated
in batch-cycle mode on various substrates with applied voltage 0.6
V at temperature $T=25^{\circ}$C. Multiple batches were conducted
for each substrate. In particular, we consider batch cycles for an
acetate fed MEC. At the beginning of the experiment, the initial concentration
of acetate was $A_{0}=956$ mg/L. During each batch cycle, current
first increases as the exoelectrogen population grows and then decreases
due to depletion of substrate. After each batch, devices were emptied
out to limit growth of methanogens. Additionally, there are only time
series measurements of current and not of hydrogen production, so
we cannot quantify the effect of hydrogenotrophic methanogens on the
hydrogen production rate. For more details about reactor construction,
see descriptions in \cite{lu_hydrogen_2009} and \cite{call_hydrogen_2008}.

\subsection{Model}

Due to the experimental considerations above, we investigate the stability
properties of a reduced version of the MEC model from Pinto et al. \cite{pinto2011}
where fermenting and hydrogenotrophic microorganisms are not considered.
We consider either a continuous-flow or a batch-cycle MEC in which
methanogens and exoelectrogens compete for a single substrate such
as acetate. In the continuous-flow case, the dilution rate, $D=F_{\text{in}}/V$,
is the influent (also effluent) flow rate divided by the reactor volume.
In the batch-cycle case, $D=0.$ The model consists of two anodic
biofilm layers. Figure \ref{fig:Cartoon_MEC} represents these microorganism
layers by two layers of spheres on the anode biofilm.  In the outer
anode biofilm layer, denoted 1, methanogenic microorganisms (depicted
as blue spheres) convert the substrate into methane and carbon dioxide.
In the inner anode biofilm layer, denoted 2, exoelectrogenic microorganisms
(depicted as green spheres) consume the substrate to produce electrons
using an intracellular mediator (or electron acceptor) while a second
compartment of methanogens competes for available substrate. In each
biofilm layer, microorganism concentration is limited by a theoretical
maximum concentration. The differential-algebraic equations representing:
substrate concentration, $S$, microorganism concentrations, $X_{m,1}$,
$X_{e}$, and $X_{m,2}$, oxidized mediator concentration, $M_{\text{ox}}$,
current, $I_{\text{MEC}}$, current density, $I_{\text{density}}$,
and internal resistance, $R_{\text{int}}$, are
\begin{align}
\frac{dS}{dt}= & D[S_{0}-S(t)]-q_{e}(t)X_{e}(t)-q_{m}(t)[X_{m,1}(t)+X_{m,2}(t)],\label{eq:acetate}\\
\frac{dX_{m,1}}{dt}= & [\mu_{m}(t)-K_{d,m}-D\alpha_{1}(t)]X_{m,1}(t),\label{eq:methanogen1}\\
\frac{dX_{e}}{dt}= & [\mu_{e}(t)-K_{d,e}-D\alpha_{2}(t)]X_{e}(t),\label{eq:exoelectrogen}\\
\frac{dX_{m,2}}{dt}= & [\mu_{m}(t)-K_{d,m}-D\alpha_{2}(t)]X_{m,2}(t),\label{eq:methanogen2}\\
\frac{dM_{\text{ox}}}{dt}= & -Y_{M}q_{e}(t)X_{e}(t)+\frac{\gamma}{VmF_{2}}I_{\text{MEC}}(t),\label{eq:oxidizedMediator}\\
I_{\text{MEC}}(t)R_{\text{int}}(t)= & E_{\text{applied}}+E_{\text{CEMF}}-\frac{RT}{mF}\ln\left(\frac{M_{\text{total}}}{M_{\text{total}}-M_{\text{ox}}(t)}\right)-\frac{RT}{\beta mF}\text{arcsinh}\left(\frac{I_{\text{MEC}}(t)}{A_{\text{sur,A}}i_{0}}\right),\label{eq:current}\\
I_{\text{density}}(t)=\frac{1000}{V} & I_{\text{MEC}}(t),\label{eq:currentDensity}\\
R_{\text{int}}(t)= & R_{\text{min}}+(R_{\text{max}}-R_{\text{min}})e^{-K_{R}X_{e}(t)},\label{eq:internalResistance}
\nonumber 
\end{align}

\begin{table}
\centering%
\scalebox{0.85}{%
\begin{tabular*}{1.1\textwidth}{|c|c|c|c|c|}
\hline 
Parameter & Description & Value & Units & Source \\
\hline 
\hline 
$D$ & dilution rate & 0 & 1 / day & experiment \\
\hline 
$S_{0}$ & influent substrate concentration & 956 & mg-$S$ / L & experiment \\
\hline 
$E_{\text{applied}}$ & applied voltage & 0.6 & Volts & experiment \\
\hline 
$A_{\text{sur,A}}$ & anode surface area & 0.4 & m\textasciicircum{}2 & experiment \\
\hline 
$i_{0}$ & equilibrium exchange current density & 1 & ampere / m\textasciicircum{}2 & \cite{pinto2011} \\
\hline 
$V$ & MEC liquid volume & 0.09 & L & experiment \\
\hline 
$T$ & MEC temperature & 298.15 & K & experiment\\
\hline 
$P$ & MEC pressure & 1 & atm & experiment\\
\hline 
$F$ & Faraday's constant (in seconds) & 96485 & ampere sec / mol-$e^{-}$ & constant\\
\hline 
$R$ & ideal gas constant (in joules) & 8.3145 & J / mol / K & constant\\
\hline 
$E_{\text{CEMF}}$ & counter-electromotive force & -0.34 & Volts & \cite{pinto2011}\\
\hline 
$\mu_{\text{max},e}$ & max. growth rate of exoelectrogen & 2.43 & 1 / day & fit\\
\hline 
$\mu_{\text{max},m}$ & max. growth rate of methanogen & 0.3 & 1 / day & assumed\\
\hline 
$q_{\text{max},e}$ & max. consumption rate, exoelectrogen & 4.82 & mg-$S$ / mg-$X_{e}$ / day & fit\\
\hline 
$q_{\text{max},m}$ & max. consumption rate, methanogen & 4 & mg-$S$ / mg-$X_{m}$ / day & assumed\\
\hline 
$K_{S,e}$ & half rate constant, exoelectrogen & 800 & mg-$S$ / L & assumed\\
\hline 
$K_{S,m}$ & half rate constant, methanogen & 810 & mg-$S$ / L & assumed \\
\hline 
$K_{M}$ & half rate constant of mediator & 0.2$M_{\text{total}}$  & mg-$M$ / L & \cite{pinto2011}\\
\hline 
$K_{d,e}$ & decay rate for exoelectrogens & 0.04 & 1 / day & \cite{pinto2011}\\
\hline 
$K_{d,m}$ & decay rate for methanogens & 0.002 & 1 / day & \cite{pinto2011}\\
\hline 
$K_{X}$ & curve steepness for biofilm retention & 0.04 & - & \cite{pinto2011}\\
\hline 
$X_{\text{max,1}}$ & max. concentration in biofilm 1 & 900 & mg-$X$ / L & \cite{pinto2011}\\
\hline 
$X_{\text{max,2}}$ & max. concentration in biofilm 2 & 512.5 & mg-$X$ / L & \cite{pinto2011}\\
\hline 
$Y_{M}$ & oxidized mediator yield & 40.7 & mg-$M$ / mg-$S$ & fit\\
\hline 
$R_{\text{min}}$ & minimum internal resistance & 25 & Ohms & \cite{pinto2011}\\
\hline 
$R_{\text{max}}$ & maximum internal resistance & 2000 & Ohms & \cite{pinto2011}\\
\hline 
$K_{R}$ & curve steepness of internal resistance & 0.06 & L / mg-$X_{e}$ & assumed\\
\hline 
$M_{\text{total}}$ & max. mediator concentration & 0.05$X_{\text{max},2}$ & mg-$M$ / mg-$X_{e}$ & \cite{pinto2011}\\
\hline 
$\gamma$ & assumed mediator molar mass & 663400 & mg-$M$ / mol-$M$ & \cite{pinto2011}\\
\hline 
$m$ & electrons transferred per mol mediator & 2 & mol-$e^{-}$ / mol-$M$ & \cite{pinto2011}\\
\hline 
$\beta$ & reduction \& oxidation transfer coefficients & 0.5 & - & \cite{pinto2011}\\
\hline
\end{tabular*}}\caption{Description of the parameters in equations (\ref{eq:acetate})-(\ref{eq:biofilm2}).
Most parameters come from physical constants, experimental measurements,
or data from \cite{pinto2011}. $\mu_{\text{max},m}$, $q_{\text{max},m}$,
$K_{S,e}$, $K_{S,m}$, and $K_{R}$ were not identifiable. Values
for these were assumed based on sensitivities and preliminary fits.
$\mu_{\text{max},e}$, $q_{\text{max},e}$, and $Y_{M}$ were fit
to the data.}\label{tab:parameter_description}
\end{table}
Table \ref{tab:parameter_description} contains a description of
the model parameters, including units. The growth rates and consumption
rates of the methanogenic and exoelectrogenic microorganisms can be
defined using Monod kinetics. Our equations differ from \cite{pinto2011}
in that we use oxidized mediator concentration, instead of the fraction
of oxidized mediator per exoelectrogen. We choose concentration because
that is the natural unit for the Monod terms.
The growth and consumption rates are

\begin{align}
\mu_{e}(t) & =\mu_{\text{max},e}\left(\frac{S(t)}{K_{S,e}+S(t)}\right)\left(\frac{M_{\text{ox}}(t)}{K_{M}+M_{\text{ox}}(t)}\right),\label{eq:mu_e}\\
\mu_{m}(t) & =\mu_{\text{max},m}\left(\frac{S(t)}{K_{S,m}+S(t)}\right),\label{eq:mu_m}\\
q_{e}(t) & =q_{\text{max},e}\left(\frac{S(t)}{K_{S,e}+S(t)}\right)\left(\frac{M_{\text{ox}}(t)}{K_{M}+M_{\text{ox}}(t)}\right),\label{eq:q_e}\\
q_{m}(t) & =q_{\text{max},m}\left(\frac{S(t)}{K_{S,m}+S(t)}\right).\label{eq:q_m}
\end{align}

The microorganism concentration in each biofilm layer is limited by
a maximum concentration parameter. The dimensionless biofilm retention
functions are continuous functions which provide a mechanism for decreasing
the rate of change in microorganism concentration in a continuous-flow
MEC when the concentration exceeds the maximum in a given layer. In
their continuous-flow MEC model, Pinto et al. used piecewise
biofilm retention functions which are zero when biomass is below the
maximum and nonzero when biomass exceeds the maximum in each layer \cite{pinto2011}.
However, we use the hyperbolic tangent formulation from their previous
MFC model \cite{pinto_two-population_2010}. We do this because
continuous equations result in a more robust numerical scheme. The
biofilm retention functions for biofilms 1 and 2 are
\begin{align}
\alpha_{1}(t) & =\frac{1}{2}\left(1+\text{Tanh}[K_{x}(X_{m,1}(t)-X_{\text{max},1})]\right),\label{eq:biofilm1}\\
\alpha_{2}(t) & =\frac{1}{2}\left(1+\text{Tanh}[K_{x}(X_{e}(t)+X_{m,2}(t)-X_{\text{max},2})]\right).\label{eq:biofilm2}
\end{align}

To derive the equation for MEC current, we follow Pinto et al.~in
using the following electrochemical balance equation:
\begin{equation}
E_{\text{applied}}+E_{\text{CEMF}}=\eta_{\text{ohm}}+\eta_{\text{act,A}}+\eta_{\text{act,C}}+\eta_{\text{conc,A}}+\eta_{\text{conc,C}}\label{eq:electrochemicalBalance}
\end{equation}
where $E_{\text{applied}}$ is the applied voltage; $E_{\text{CEMF}}$
is the counter-electromotive force; $\eta_{\text{ohm}}$ is the ohmic
loss; $\eta_{\text{act,A}}$ and $\eta_{\text{act,C}}$ are the activation
losses at the anode and cathode, respectively; and $\eta_{\text{conc,A}}$
and $\eta_{\text{conc,C}}$ are the concentration losses at the anode
and cathode, respectively. Note that activation and concentration
losses apply at both the anode and the cathode. Pinto et al. neglect
concentration losses at the cathode due to the assumption that hydrogen
molecules diffuse away from the cathode rapidly. The authors also
assume that activation losses can be neglected at the anode since
the MEC operates at high overpotential at the cathode. This allows
us to deal with only two nonlinear terms instead of four. Ohmic losses
can be calculated from Ohm\textquoteright s Law: $\eta_{\text{ohm}}=I_{\text{MEC}}R_{\text{int}}$.
Following Marcus et al. \cite{kato_marcus_conduction-based_2007}, the authors
write concentration losses at the anode using the Nernst equation
with the assumption that the reference reduced mediator concentration
(or standard anodic electron acceptor concentration) is equal to the
total intracellular mediator concentration \cite{pinto_two-population_2010,pinto2011}.
Then
\begin{equation}
\eta_{\text{conc, A}}=\frac{R_{1}T}{mF_{1}}\ln\left(\frac{M_{\text{total}}}{M_{\text{red}}}\right).\label{eq:concentrationLossAnode}
\end{equation}
Pinto et al. also calculate activation losses at the cathode
using the Butler-Volmer equation which relates potential to current
at an electrode \cite{pinto2011}. We use standard simplifying assumptions that the
reaction occurs in one step and that the symmetry coefficient (or
the fraction of activation loss that affects the rate of electrochemical
transformation) is $\beta=0.5$. With these assumptions we can write
\begin{equation}
\eta_{\text{act, C}}=\frac{R_{1}T}{\beta mF_{1}}\text{arcsinh}\left(\frac{I_{\text{MEC}}}{A_{\text{sur,A}}i_{0}}\right).\label{eq:activationLossCathode}
\end{equation}
For more information, see the explanation in \cite{noren_clarifying_2005}
about approximations to the Butler-Volmer equation. Equations (\ref{eq:electrochemicalBalance})-(\ref{eq:activationLossCathode})
combine to give current implicitly from the nonlinear function in
equation (\ref{eq:current}).

\subsection{Parameter Fitting}

Parameter fitting was performed using the trust-region-reflective
algorithm as implemented in MATLAB's nonlinear least-squares solver
\textit{lsqnonlin}. The solution sensitivities discussed in Section
\ref{sec:Sensitivity-Equations} were used to specify an objective
gradient and to determine which parameters are identifiable. The fitted
parameters for each batch of the acetate fed MEC are shown in Table
\ref{tab:parameter_description}. Some parameters are physical constants,
others are known from the batch-cycle experiment described above,
and some were taken from values in \cite{pinto2011}. For the other
parameters, we determined that the maximum exoelectrogen growth rate,
$\mu_{\text{max},e}$, maximum exoelectrogen consumption rate, $q_{\text{max},e}$,
and mediator yield, $Y_{M}$, are identifiable. We estimated these
parameters starting with values reported in the supplementary table
for \cite{pinto2011}. However, the values from the supplementary
table for the half rate constants, $K_{S,e}$ and $K_{S,m}$, and
resistance steepness, $K_{R}$, did not provide a good fit to the
data. In addition, these three parameters were unidentifiable. Therefore,
their values were assumed based on sensitivities and preliminary fits
to provide a reasonable initial guess for the nonlinear least squares
algorithm.

\section{Sensitivity Equations\label{sec:Sensitivity-Equations}}

To perform differential-algebraic sensitivity analysis, we solve the
sensitivity equations for the model. In general, a DAE with parameters
$\mathbf{p}$ can be written as
\[
\mathbf{F}(t,\mathbf{y},\mathbf{y}',\mathbf{p})=0,\quad\mathbf{y}(t_{0})=\mathbf{y}_{0},\quad\dot{\mathbf{y}}(t_{0})=\dot{\mathbf{y}}_{0}.
\]
Let $\mathbf{s}_{i}(t)$ denote the solution sensitivity with respect
to the parameter $p_{i}$. That is, $\mathbf{s}_{i}(t)=\frac{\partial\mathbf{y}(t)}{\partial p_{i}}$.
Then the sensitivity equations, with respect to parameter $p_{i}$,
can be written as
\begin{align*}
\frac{\partial\mathbf{F}}{\partial\mathbf{y}}\mathbf{s}_{i}+\frac{\partial\mathbf{F}}{\partial\dot{\mathbf{y}}}\dot{\mathbf{s}}_{i}+\frac{\partial\mathbf{F}}{\partial p_{i}}=0,\\
\mathbf{s}_{i}(t_{0})=\frac{\partial\mathbf{y}_{0}(\mathbf{p})}{\partial p_{i}},\quad\dot{\mathbf{s}}_{i}(t_{0})=\frac{\partial\dot{\mathbf{y}}_{0}(\mathbf{p})}{\partial p_{i}}.
\end{align*}
The precise form of these equations will depend on the parameter of
interest. We are interested in the solution sensitivity of current
density with respect to parameter $p_{i}$, that is $\frac{\partial}{\partial p_{i}}I_{\text{density}}(t,\mathbf{p}),$
which is the final component of $\mathbf{s}_{i}(t)$. Since the parameters
have different units and magnitudes, the solution sensitivities are
not directly comparable. To remedy this, we will look at the semi-relative
sensitivity of current density, $p_{i}\frac{\partial}{\partial p_{i}}I_{\text{density}}(t,\mathbf{p})$.
The semi-relative sensitivity of current density to any parameter
has units of ampere/m\textasciicircum{}3. This allows us to compare
how the current density changes density changes with respect to changes
in parameters. 

The model and sensitivity system were solved simultaneously using
a variable-order, variable-coefficient backward differentiation formula
in fixed-leading coefficient form \cite{brenan_numerical_1995},
as implemented in the IDAS package from the SUNDIALS suite of nonlinear
and differential-algebraic equation solvers \cite{hindmarsh2005sundials}.
The initial conditions for acetate concentration, current density,
and methanogen concentration were determined by the experiments. Initial
acetate concentration was fixed at 956 mg/L and initial current density
was measured as 45.8 amp/m\textasciicircum{}3 for the data used in
the fitting procedure. Initial methanogen concentrations were assumed
to be 10 mg/L because the experiment was designed so that methanogen
concentration would be significantly smaller than exoelectrogen concentration.
Initial oxidized mediator was chosen to be 25.6 mg/L to provide a
smooth solution curve. Finally, we solved the algebraic constraint
for exoelectrogen concentration to provide consistent initial conditions
for the DAE. This was done using the trust-region-dogleg algorithm
as implemented in MATLAB's nonlinear equation solver \textit{fsolve}.
We also verified that the initial conditions were consistent using
Newton iteration paired with a global line search strategy as implemented
in the IDACalcIC routine in IDAS. Newton corrections made use of a
dense linear solver and a user specified Jacobian. According to both
methods, the initial guess $X_{e}(0)=250$ mg/L satisfied the constraint.

\begin{figure}
\centering{}\includegraphics[width=0.75\textwidth]{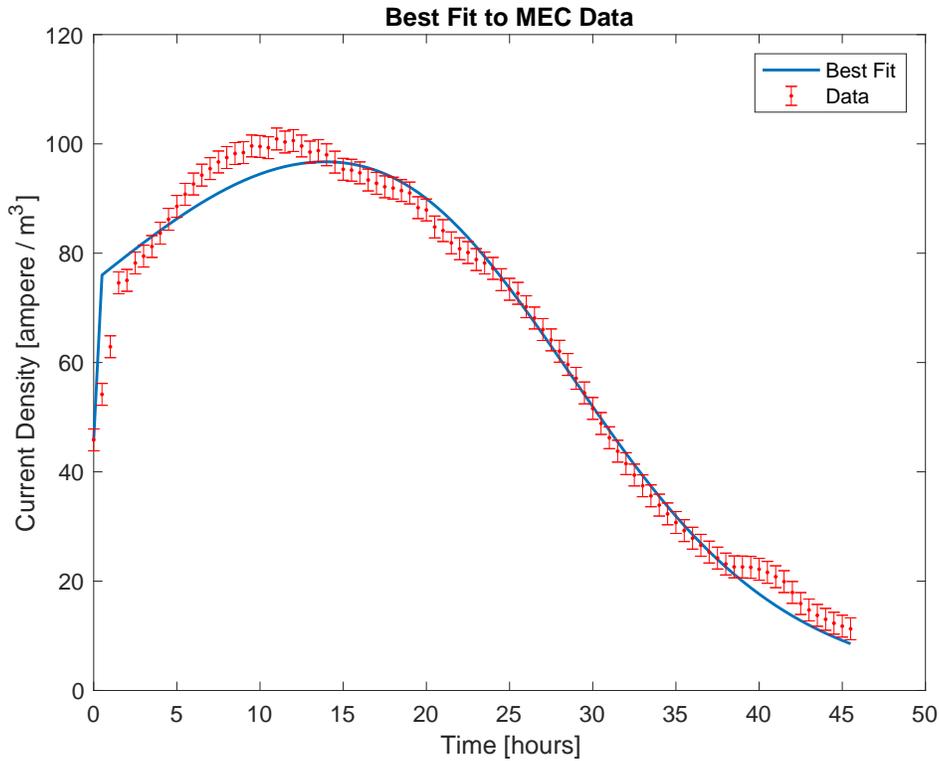}\caption{Comparison of best fit current density with data from a batch of the
acetate fed MEC. The best fit curve for current density (in ampere/m\textasciicircum{}3) is shown as a solid blue line. Experimental
data is depicted as a series of red circles. Error bars show the measurement
error, $\pm2$ ampere/m\textasciicircum{}3. \label{fig:current_density_fit-2}}
\end{figure}

\begin{figure}
\centering{}\includegraphics[width=0.75\textwidth]{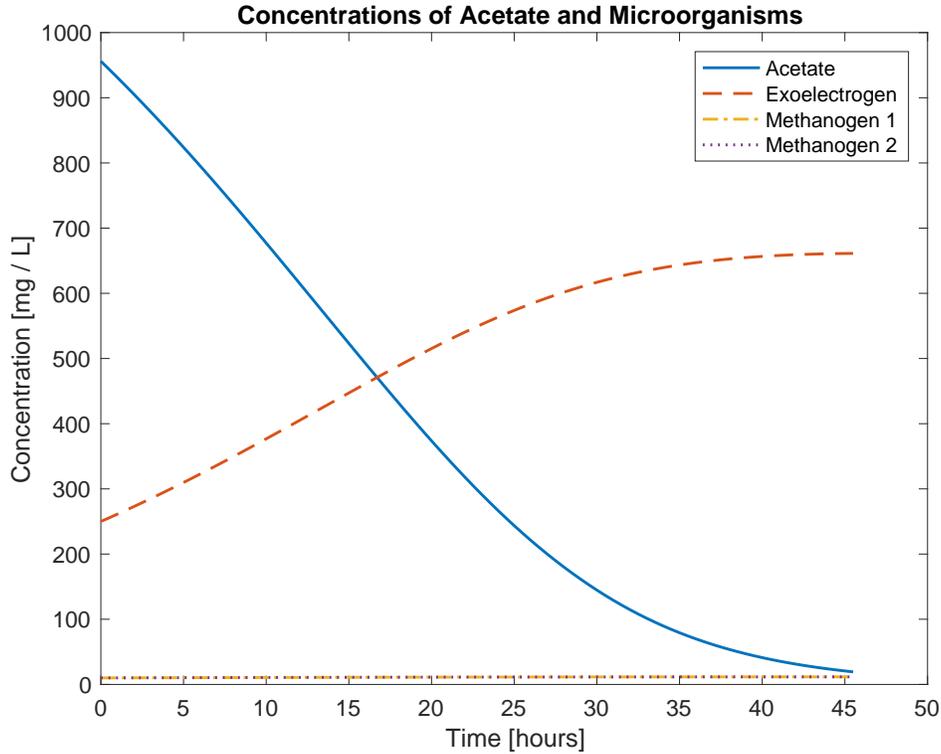}\caption{Best fit concentrations of acetate, exoelectrogens, and methanogens
in mg/L. Acetate is represented by a solid blue line. Exoelectrogen
concentration is shown as a dashed red line. Methanogen concentrations
in biofilms 1 and 2 are depicted by a dashed-dot yellow line and a
dotted purple line, respectively. \label{fig:concentrations-1}}
\end{figure}

\begin{figure}
\centering{}\includegraphics[width=0.75\textwidth]{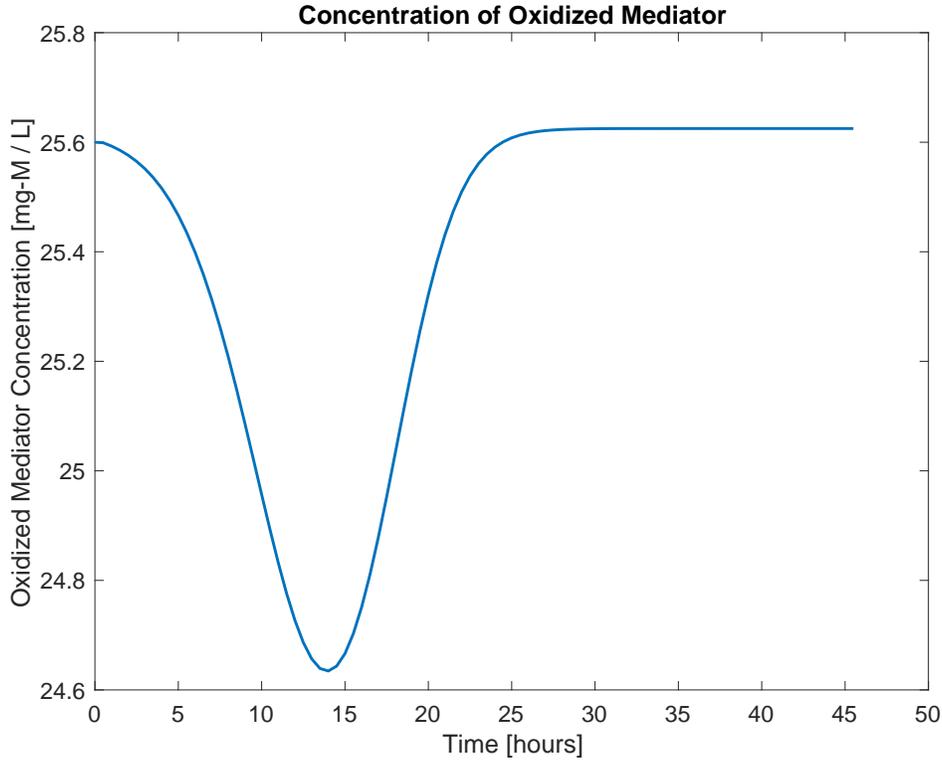}\caption{Solution of oxidized mediator concentration (out of $M_{\text{total}}=25.625$
mg/L). The initial conditions were chosen by solving the algebraic
equation for the initial exoelectrogen concentration, $X_{e}(0)$,
given initial conditions for oxidized mediator concentration and current
density. $I_{\text{density}}(0)$ was experimentally determined. $M_{\text{ox}}(0)=25.6$
was chosen to provide a smooth curve for the oxidized mediator concentration.
\label{fig:oxidized_mediator-1}}
\end{figure}

Figure \ref{fig:current_density_fit-2} shows the best fit current
density, $I_{\text{density}}(t)$, compared to experimental data for
one batch with an acetate fed MEC. Other batches provided similar
results. The simulations consistently underestimate or overestimate
the data for several hours at a time, but provide a reasonable approximation
to the current density over a few days. In particular, the simulation
underestimates the data from 5 to 14 hours and again after 39 hours.
During these time periods, there is some feature of the model that
is not captured by the mathematics as well as we would like. Figure
\ref{fig:concentrations-1} shows the solution for acetate, exoelectrogen,
and methanogen concentrations and Figure \ref{fig:oxidized_mediator-1}
shows the solution for oxidized mediator concentration. We do not have
data for state variables besides current density, so these figures show only simulations. 

\begin{figure}
\centering{}\includegraphics[width=0.75\textwidth]{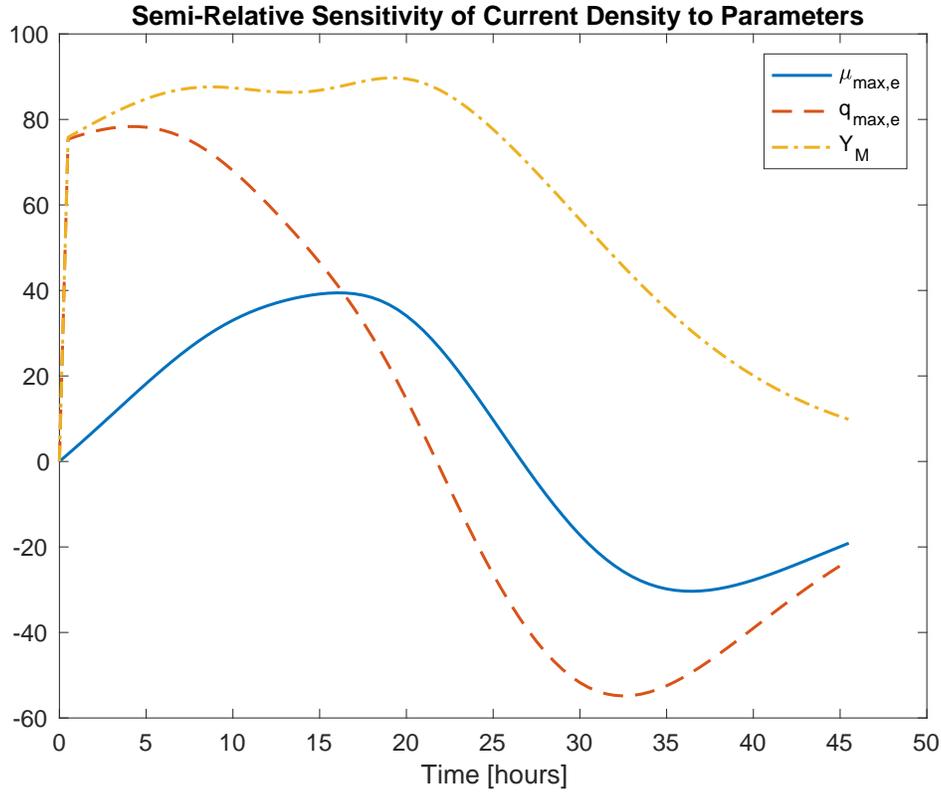}\caption{Semi-relative sensitivity of current density with respect to fitted
parameters. Sensitivity to $\mu_{\text{max},e}$ is shown as a solid
blue line, sensitivity to $q_{\text{max},e}$ is depicted as a dashed
red line, and sensitivity to $Y_{M}$ is presented as a dashed-dot
yellow line. Increasing either $q_{\text{max},e}$ or $Y_{M}$ will
significantly increase current density during the first few hours
of the experiment. This effect tapers off over time for $Y_{M}$.
For $q_{\text{max},e}$, the effect is eventually reversed. Between
$t=22$ and 45 hours, increasing $q_{\text{max},e}$ will decrease
current density. Increasing $\mu_{\text{max},e}$ will increase current
density before $t=27$ and decrease current density afterwards, although
the effects are less pronounced. \label{fig:semi-relative_sensitivity-1}}
\end{figure}

\begin{figure}
\begin{centering}
\includegraphics[width=0.75\textwidth]{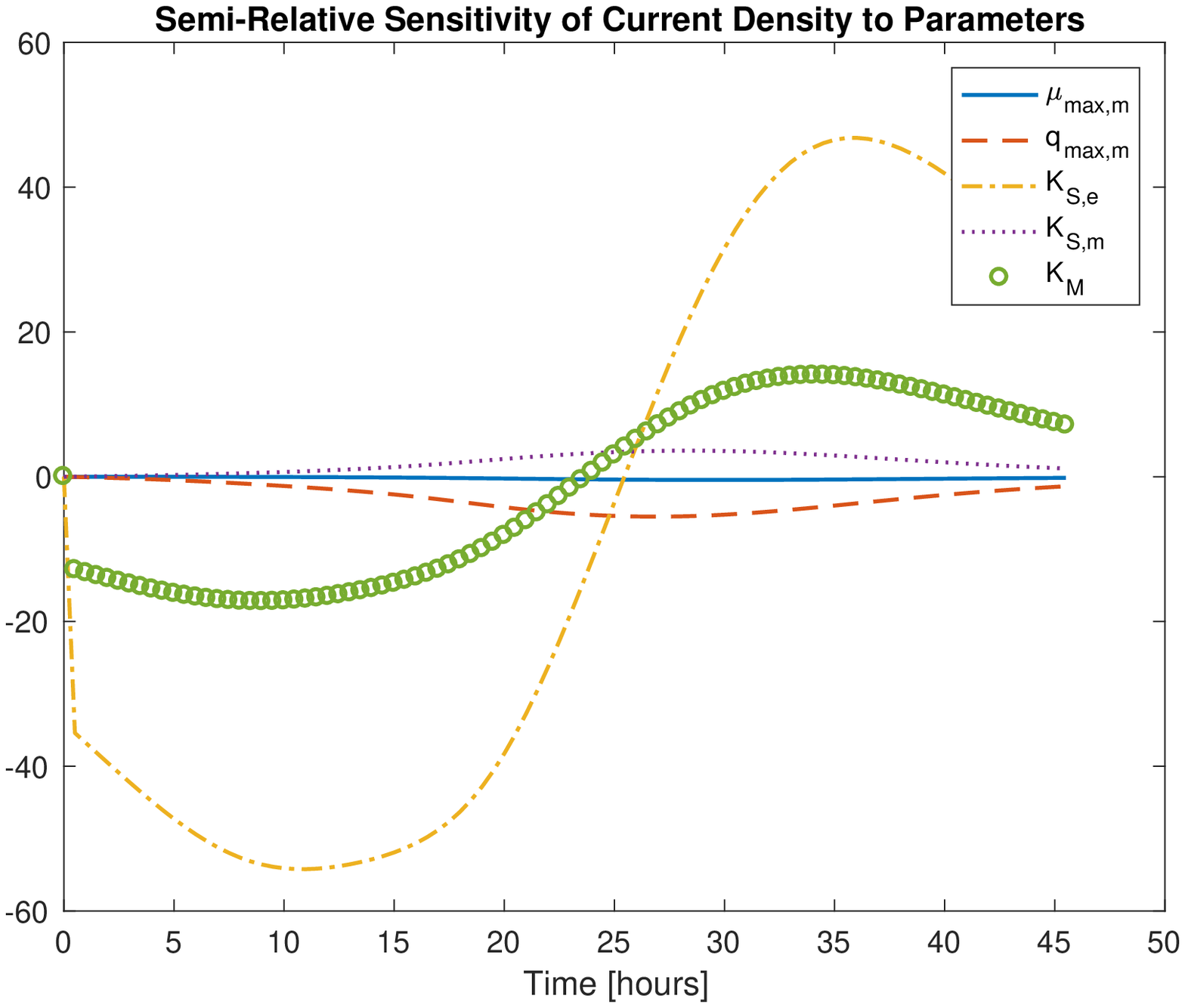}
\par\end{centering}
\caption{Semi-relative sensitivity of current density with respect to methanogen
growth and consumption parameters and the three half rate constants.
These parameters could not be identified. Sensitivity $\mu_{\text{max},m}$
is shown as a solid blue line, sensitivity to $q_{\text{max},m}$
is depicted as a dashed red line, sensitivity to $K_{S,e}$ is pictured
as a dashed-dot yellow line, sensitivity to $K_{S,m}$ is pictured
as a dotted purple line, and sensitivity to $K_{M}$ is pictured as
a series of green circles. $\mu_{\text{max},m}$, $q_{\text{max},m}$,
and $K_{S,m}$ do not have a significant influence on current density
because methanogen concentration is relatively small in the experiment.
Increasing $K_{S,e}$ would significantly decrease the current density
before $t=25$ and would increase current density later. Modifying
$K_{M}$ has similar but more moderate effects. \label{fig:semi-relative-sensitivity-2}}
\end{figure}

The main sensitivity analysis results are displayed in Figures \ref{fig:semi-relative_sensitivity-1}
and \ref{fig:semi-relative-sensitivity-2}. Each curve represents
the semi-relative sensitivity of current density with respect to a
particular parameter. We interpret a curve as the influence of the
corresponding parameter on the current density at various times. $\mu_{\text{max},e}$,
$q_{\text{max},e}$, $K_{S,e}$, $K_{M}$, and $Y_{M}$ have the strongest
influence on current density throughout the experiment. These are
the parameters related to exoelectrogen growth and consumption. Only
$\mu_{\text{max},e}$, $q_{\text{max},e}$, and $Y_{M}$ were identifiable
for parameter fitting. Since the derivative with respect to the consumption
rate, $q_{\text{max},e}$, and the oxidized mediator yield, $Y_{M}$,
are large during the first few hours, increasing either of these parameters
would significantly increase the current density during this time.
This effect continues for $Y_{M}$, but tapers off over the course
of the experiment. The effect of increasing $q_{\text{max},e}$ is
not obvious beforehand. However, Figure \ref{fig:semi-relative_sensitivity-1}
shows that increasing $q_{\text{max},e}$ would cause the current
density to increase during the first 25 hours but would cause it to
decrease after that. As one might expect, increasing $\mu_{\text{max},e}$
would increase the exoelectrogen growth rate and lead to increases
in current density, while increasing $K_{S,e}$ would decrease the
the exoelectrogen growth rate and lead to decreases in current density,
at least during the first day. During the second day, increasing these
parameters would have the opposite effects since they promote either
depletion or retention of substrate during the first day, The effects
of $K_{S,e}$ mirror those of $\mu_{\text{max},e}$ and $q_{\text{max},e}$,
due to the structure of the Monod kinetics equations. Perhaps unsurprisingly,
the methanogen growth rate parameters, $\mu_{\text{max},m}$, $q_{\text{max},m}$,
and $K_{S,m}$, do not have much influence on current during the first
45 hours. This is likely related to the low concentration of methanogens
in this experiment, Methanogen populations remain small during the
experiment, but could become dominant if the simulations were run
for longer periods of time. This competition is an interesting feature
of the continuous-flow model and its implications are discussed further
in Section \ref{sec:Bifurcation-Analysis}.


\section{Bifurcation Analysis in Dilution Rate\label{sec:Bifurcation-Analysis}}

Before presenting the bifurcation results, we provide a brief review
of stability of equilibria and bifurcations in DAEs. The model in
equations (\ref{eq:acetate}) - (\ref{eq:current}) can be represented
as a regular, semi-explicit, index-1 DAE
\begin{align}
y'= & \,f(y,z),\label{eq:semi-explicit-1}\\
0= & \,g(y,z),\nonumber 
\end{align}
where $f:\mathbb{R}^{n+m}\rightarrow\mathbb{R}^{n}$ and $g:\mathbb{R}^{n+m}\rightarrow\mathbb{R}^{m}$.
We will call a point \emph{regular} for the semi-explicit DAE (\ref{eq:semi-explicit-1})
if and only if $g(y,z)=0$ and $g_{z}(y,z)$ defines a nonsingular
matrix \cite{riaza_stability_2002,riaza_2008}. At a regular point,
we can differentiate the algebraic part of (\ref{eq:semi-explicit-1}),
as long as $g\in C^{1}$, to obtain what is known as the underlying
ODE
\begin{align*}
y' & =f(y,z)\\
z' & =-g_{z}^{-1}(y,z)g_{y}(y,z)f(y,z).
\end{align*}
The fact that only one differentiation is required at a regular point
defines (\ref{eq:semi-explicit-1}) as an index-1 DAE. Alternatively,
regularity at $(y,z)$ allows us to use the implicit function theorem
to conclude the existence of a map $z=\phi(y)$ near $(y,z)$ and
to describe behavior on the differential solution manifold $g(y,z)=0$
by 
\begin{equation}
y'=f(y,\phi(y)).\label{eq:reduced}
\end{equation}
 This is known as the reduced ODE \cite{riaza_stability_2002,beardmore_stability_1998}.

Equilibrium points of (\ref{eq:semi-explicit-1}) satisfy both $f(y,z)=0$
and $g(y,z)=0$. After defining $\hat{f}(y)=f(y,\phi(y))$, one can
use the implicit function theorem and derivatives of (\ref{eq:semi-explicit-1})
to show that linear stability of equation (\ref{eq:reduced}) at regular
equilibria is determined by the Schur complement of the lower right
block of the Jacobian matrix of $[f,g]^{T}$ \cite{venkatasubramanian_local_1995,beardmore_stability_1998,reich_local_1995}.
In particular,
\begin{equation}
\hat{f}_{y}=f_{y}-f_{z}g_{z}^{-1}g_{y}.\label{eq:Schur}
\end{equation}
 Therefore a regular equilibrium of (\ref{eq:semi-explicit-1}) is
asymptotically stable if all the eigenvalues of $f_{y}-f_{z}g_{z}^{-1}g_{y}$
have negative real part. Computing $f_{y}-f_{z}g_{z}^{-1}g_{y}$ is
problematic at \emph{singular} points where $g_{z}(y,z)$ is not invertible.
However, the spectrum of this matrix is the same as the spectrum of
the matrix pencil $\{E,-J\}=\{\lambda E-J:\lambda\in\mathbb{C}\}$
given by
\[
E=\begin{bmatrix}I & 0\\
0 & 0
\end{bmatrix}\text{ and }J=\begin{bmatrix}f_{y} & f_{z}\\
g_{y} & g_{z}
\end{bmatrix}.
\]
One could also use the fact that an equilibrium of the DAE (\ref{eq:semi-explicit-1})
is asymptotically stable if the spectrum of the matrix pencil, $\sigma(\{E,-J\})=\{\lambda\in\mathbb{C}$
: $\text{det}(\lambda E-J)=0\}$, consists of elements with negative
real part \cite{riaza_stability_2002}. It is worth noting that
the matrix pencil can be used to characterize stability for a generic
DAE. For more information, see the exposition by Rabier and Rheinboldt \cite{ciarlet_1990}
or the discussion with applications to circuits by Riaza \cite{riaza_2008}.
For the MEC model, 
\[
g_{z}(y,z)=-R_{\text{int}}-\frac{RT}{\beta mFA_{\text{sur},A}i_{0}\sqrt{1+\left(\frac{I_{\text{MEC}}}{A_{\text{sur},A}i_{0}}\right)^{2}}}<0
\]
is always nonsingular. Therefore, if $f(y,z)=g(y,z)=0$, we must have
a regular equilibrium point and the equilibrium is asymptotically
stable if the eigenvalues of $f_{y}-f_{z}g_{z}^{-1}g_{y}$ all have
negative real part. 

In a batch-cycle MEC, the dilution rate, $D=F_{\text{in}}/V$, is
zero and equation (\ref{eq:acetate}) tells us that $\frac{dS}{dt}\leq0$.
In particular, $\frac{dS}{dt}=0$ only if either $S(t)=0$ or $X_{m,1}(t)=X_{e}(t)=X_{m,2}(t)=0$.
If $S(t)=0$, then for each bacteria compartment $i$, $X_{i}(t)=X_{i}(0)e^{-K_{d,i}t}$.
In either case $\lim_{t\rightarrow\infty}X_{i}(t)=0$. The remaining
system becomes 
\begin{align*}
\frac{dM_{\text{ox}}}{dt} & =\frac{\gamma}{VmF_{2}}I_{\text{MEC}}(t),\\
I_{\text{MEC}}(t)R_{\text{int}}(t) & =E_{\text{applied}}+E_{\text{CEMF}}-\frac{R_{1}T}{mF_{1}}\ln\left(\frac{M_{\text{total}}}{M_{\text{total}}-M_{\text{ox}}(t)}\right)-\frac{R_{1}T}{\beta mF_{1}}\text{arcsinh}\left(\frac{I_{\text{MEC}}(t)}{A_{\text{sur,A}}i_{0}}\right).
\end{align*}
For the parameters in Table \ref{tab:parameter_description} , the
system has a line of equilibria at 
\[
\{S,X_{m,1},X_{e},X_{m,2},M_{\text{ox}},I_{\text{MEC}}\}=\{S,0,0,0,M^{*},0\}
\]
where $M^{*}$ is approximately equal to, but less than $M_{\text{total}}$
to satisfy the constraint. For instance, with the initial conditions,
$\{S(0),X_{m,1}(0),X_{e}(0),X_{m,2}(0),M_{\text{ox}}(0),I_{MEC}(0)\}=\{956,10,250,10,5,0.006123\}$,
the system approaches
\[
\{S,X_{m,1},X_{e},X_{m,2},M_{\text{ox}},I_{\text{MEC}}\}=\{0,0,0,0,M^{*},0\}
\]
with $|M^{*}-M_{\text{total}}|=4.16\times10^{-8}$. This is a nonhyperbolic
equilibrium. The spectrum of the matrix pencil (also of the Schur
complement) has one zero eigenvalue:
\[
\sigma(\{E,-J\})=\{-5.099\times10^{8},-0.04,-0.002,-0.002,0\}.
\]
This point lies on a line of stable equilibria in the direction $\{1,0,0,0,0,0\}$
with $\frac{dS}{dt}=0$. The substrate can take any value at equilibrium
if the microorganism concentrations are zero, so any point
\[
\{S,X_{m,1},X_{e},X_{m,2},M_{\text{ox}},I_{\text{MEC}}\}=\{S^{*},0,0,0,M^{*},0\}
\]
can be a stable equilibrium point for the MEC system when $D=0$. 

\begin{figure}
\begin{centering}
\includegraphics[width=\textwidth]{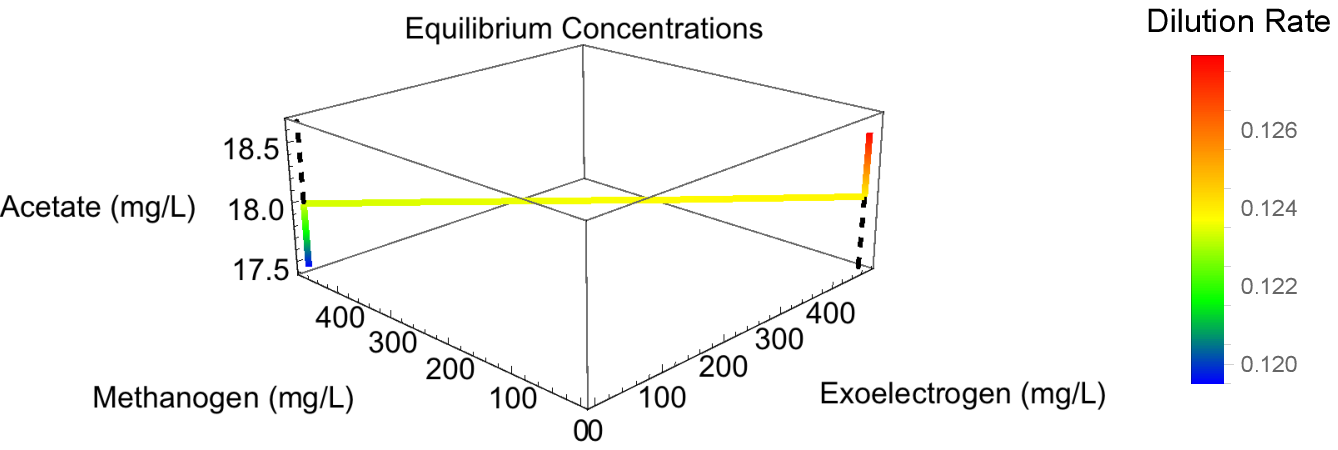}
\par\end{centering}
\caption{Stable equilibria are depicted by solid colored lines and unstable
equilibria are depicted by black dashed lines. The value of the dilution rate, $D$, increases from left to right as shown by the color bar.  As
$D$ increases, the first transcritical bifurcation moves the stable
equilibrium from a curve with no exoelectrogens to the curve of coexistence.
The second transcritical bifurcation moves the stable equilibrium
from the line of coexistence to a line with no methanogens. In other
words, increasing D moves the stable equilibrium through three regimes:
(1) competitive exclusion by methanogens, (2) coexistence, and (3)
competitive exclusion by exoelectrogens.\label{fig:transcriticalLine}}
\end{figure}

\begin{figure}
\begin{centering}
\includegraphics[width=\textwidth]{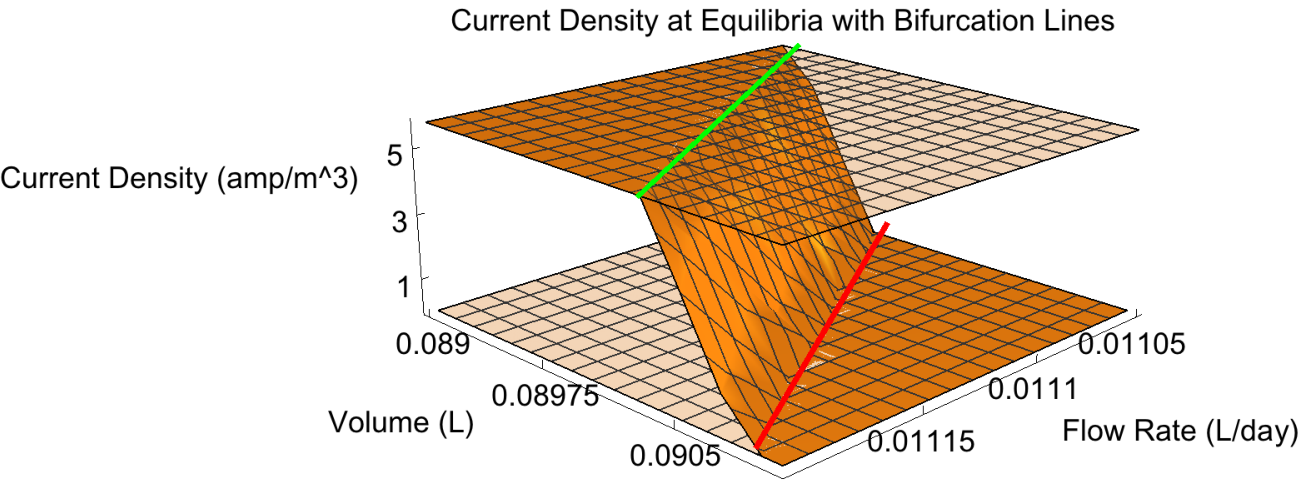}
\par\end{centering}
\caption{Stable equilibria are depicted by the solid surface and unstable equilibria
are depicted by the transparent planes. The dilution rate equals the
flow rate divided by the volume, D = Fin/V ; increasing either of
these parameters passes the stable equilibrium through three regimes,
separated by two transcritical bifurcations. Note that the flow rate
and volume axes show decreasing values for a better perspective. For
low values of D, there is competitive exclusion by methanogens, so
the equilibria lie in a plane where the current is zero. As D increases,
a transcritical bifurcation, depicted by the red line, moves the stable
equilibrium onto a plane of coexistence where current density increases
with D. As D increases further, a second transcritical bifurcation,
depicted by the green line, moves the stable equilibrium onto a plane
with competitive exclusion by exoelectrogens where the current density
at the stable equilibrium is more than 5 amp / m$^{3}$.\label{fig:transcriticalPlane}}
\end{figure}

In a continuous-flow MEC, the dilution rate is non-zero and equation
(\ref{eq:acetate}) tells us that there may be equilibria with positive
concentrations of acetate. In this scenario, we hope to find a stable
equilibrium with positive current, so that we can maintain long term
current density and hydrogen production. We expect from simple mathematical
models of chemostats that such equilibria exist for large enough flow
rates \cite{hsu_mathematical_1977,smith_theory_1995}. We demonstrate
that two transcritical bifurcations in the dilution rate parameter
cause this equilibrium to switch from competitive exclusion by methanogens
in biofilm 2 to coexistence and then to competitive exclusion by exoelectrogens.
When the dilution rate is large enough, the stable equilibrium has
positive nonzero current density of at least 5 ampere/m\textasciicircum{}3.
These results are depicted in Figures \ref{fig:transcriticalLine}
and \ref{fig:transcriticalPlane}.


Consider the parametrized DAE, 
\begin{align}
y'= & f(y,z,p),\label{eq:semi-explicit-parameterized}\\
0= & g(y,z,p).\nonumber 
\end{align}
From Sotomayor's theorem, we know that a saddle node bifurcation occurs
at a nonhyperbolic equilibrium with a geometrically simple eigenvalue
when certain nondegeneracy and transversality conditions are satisfied
\cite{sotomayor_generic_1973,meiss_differential_2017,perko_differential_2001}.
When the transversality condition is not satisfied, a transcritical
bifurcation may occur. For DAEs, we do not necessarily know the Jacobian
of the reduced ODE (\ref{eq:reduced}) near an equilibrium $(y^{*},z^{*},p^{*})$.
However, we can apply the bifurcation conditions to the Schur complement,
$f_{y}-f_{z}g_{z}^{-1}g_{y}$, as discussed in equation (\ref{eq:Schur}).
To be precise, for the parametrized DAE (\ref{eq:semi-explicit-parameterized}),
suppose that $f_{y}-f_{z}g_{z}^{-1}g_{y}$ has one geometrically simple
zero eigenvalue that is the only eigenvalue on the imaginary axis.
Furthermore, suppose that this eigenvalue has right eigenvector $v$
and left eigenvector $w$. Then a transcritical bifurcation occurs
in the parameter $p$ at $(y^{*},z^{*},p^{*})$ if the following three
conditions are satisfied. 

\begin{condition}\label{transversality-1}   
$w^{T}(f_{p}-f_{y}g_{y}^{-1}g_{p})=0$. 
\end{condition}
\begin{condition}\label{transversality-2}   
$w^{T}[(f_{p}-f_{y}g_{y}^{-1}g_{p})_{y}v]\neq0$. 
\end{condition}
\begin{condition}\label{nondegeneracy}   
$w^{T}[(f_{y}-f_{z}g_{z}^{-1}g_{y})_{y}(v,v)]\neq0$. 
\end{condition}

With the parameters fit to the MEC data, $\{\mu_{\text{max},e},q_{\text{max},e},Y_{M}\}=\{2.43,4.82,40.7\}$,
a transcritical bifurcation occurs in the dilution rate, at $D=0.1233388$.
Below this parameter value, the stable equilibrium contains only methanogens
in biofilm 2. However, at this parameter value, the equilibrium point
\[
\{S,X_{m,1},X_{e},X_{m,2},M_{\text{ox}},I_{\text{MEC}}\}=\{17.996,859.14,0,471.64,25.625,0\}
\]
has eigenvalues 
\[
\{-5.099\times10^{8},-6.406,-0.3037,-0.1666,3.036\times10^{-17}\}
\]
where the smallest eigenvalue has left eigenvector $w=<0,0,1,0,0>^{T}$
and right eigenvector
\[
v=<-0.0001481,-0.0001507,0.7070,-0.7072,-4.982\times10^{-9}>.
\]
The four conditions for a transcritical bifurcation are satisfied
since the smallest eigenvalue is zero up to machine precision and
\begin{align*}
w^{T}(f_{p}-f_{y}g_{y}^{-1}g_{p}) & =0,\\
w^{T}[(f_{p}-f_{y}g_{y}^{-1}g_{p})_{y}v] & =0.9161,\\
w^{T}[(f_{y}-f_{z}g_{z}^{-1}g_{y})_{y}(v,v)] & =-0.0003499.
\end{align*}
Beyond this dilution rate, the stable equilibrium briefly lies on
a curve of coexistence between methanogens and exoelectrogens in biofilm
2. However, another transcritical bifurcation at $D=0.1240194$ switches
stability to a third curve with only exoelectrogens in biofilm 2.
Figure \ref{fig:transcriticalLine} shows the stable equilibrium concentrations
moving between three lines in (Exoelectrogen, Methanogen, Acetate)
space as $D$ increases. Unstable equilibria are depicted as dashed
lines. Figure \ref{fig:transcriticalPlane} provides another perspective,
showing that the stable equilibrium moves between three planes in
(Flow Rate, Volume, Current Density) space as $D$ increases. In both
cases, the stable equilibria move through three regimes as $D$ increases:
(1) competitive exclusion by methanogens, (2) coexistence, and (3)
competitive exclusion by exoelectrogens. For these 90 mL MECs, the
dilution rate $D=F_{\text{in}}/V=0.124$ corresponds to a flow rate
of about $F_{\text{in}}=11$ mL / day. However, the location of the
bifurcations depends very much on the parameters used.

\section{Discussion\label{sec:Discussion}}

The differential-algebraic sensitivity analysis in Section \ref{sec:Sensitivity-Equations}
reveals how perturbations of a single parameter influence the current
density at particular times. During the first day of the experiment,
the four parameters related to exoelectrogen growth, $\mu_{\text{max},e}$,
$q_{\text{max},e}$, $K_{S,e}$, and $Y_{M}$ had the most influence
on current. Increasing the oxidized mediator yield, $Y_{M}$, would
significantly increase current throughout the experiment. In contrast,
increasing the maximum exoelectrogen consumption rate, $q_{\text{max},e}$,
would significantly increase the current during the first 22 hours,
but decrease it later; changing the exoelectrogen half rate constant,
$K_{S,e}$ would have the opposite effects of changing $q_{\text{max},e}$.
Increasing the maximum exoelectrogen growth rate, $\mu_{\text{max},e}$,
would moderately increase current during the first day, but decrease
it later; increasing the mediator half rate constant, $K_{M}$, would
have opposite effects of increasing $\mu_{\text{max},e}$ on a smaller
scale. These parameter effects show that the microbial growth was closely
correlated with substrate availability. The parameters related to
methanogen growth, $\mu_{\text{max},m}$, $q_{\text{max},m}$, and
$K_{S,m}$, had very little impact on current during experiment. This
is likely due to the small methanogen population in the MECs, which
is also consistent with previous findings that exoelectrogens have
more affinities than acetolastic methanogens \cite{lu_microbial_2016}.
These results suggest that increasing $\mu_{\text{max},e}$, $q_{\text{max},e}$,
and $Y_{M}$ or decreasing $K_{S,e}$ would result in higher peak
current density. Although the methanogen populations are small in
this experiment, the results also suggest that decreasing $q_{\text{max},m}$
or increasing $K_{S,m}$ could lead to higher peak current density
when the methanogen population is large. 

The bifurcation analysis in Section \ref{sec:Bifurcation-Analysis}
shows that the dilution rate in a continuous flow MEC must be chosen
carefully to ensure that the system approaches a stable equilibrium
with positive current density. For the fitted parameters, the system
exhibits two transcritical bifurcations. As dilution rate increases,
the first bifurcation moves the stable equilibrium from a curve with
competitive exclusion by methanogens to a curve with coexistence of
exoelectrogens and methanogens. The second bifurcation switches stability
from the curve with coexistence to a curve with competitive exclusion
by exoelectrogens. This result is depicted in Figures \ref{fig:transcriticalLine}
and \ref{fig:transcriticalPlane}. Only at high enough dilution rates
will exoelectrogens be able to dominate and provide positive current
density. This indicates that using the appropriate dilution rate can
be an effective approach for maintaining a vibrant exoelectrogenic
microbial community and maintaining stable system performance. For
the 90 mL MECs considered, these bifurcations occur at dilution rates
$D=F_{\text{in}}/V$ of about 0.1233 and 0.1240, corresponding to
flow rates just above $F_{\text{in}}=11$ mL / day. However, the location
of these bifurcations is determined by the microbial growth parameters.
Although not discussed in the results, modifying parameters besides
dilution rate may increase or decrease the current at the stable equilibrium. 

\section{Conclusion\label{sec:Conclusion}}

The microbial growth and consumption parameters were not measured in the
batch-cycle MEC experiment, but the results suggest that they should
be considered carefully in MEC studies. These parameters, namely $\mu_{\text{max},e}$,
$q_{\text{max},e}$, $K_{S,e}$, $\mu_{\text{max},m}$, $q_{\text{max},m}$,
$K_{S,m}$, and $Y_{M}$, should be estimated or measured if possible,
in both batch-cycle and continuous-flow MECs. In the former case,
the parameter values provide information about the peak current density
during each batch. In the latter case, the parameters can be used
to guide design before long term operation of a continuous-flow MEC.
Additionally, we have not considered the effects of different substrates
and different combinations of microorganisms. If we consider another
substrate, such as ethanol, the best fit values of the parameters
will be different and the sensitivity analysis methods could provide
additional insight. A further complication is that microbial community
structure varies when different substrates are used. That is, the
percentages of bacteria species that are present will depend on the
substrate that is used. This means that estimates for growth and consumption
parameters from one MEC study may not be reliable in another. These
factors are not taken into account by this simple model, but a similar
modeling approach can be used as the degradation process is comparable.
We plan to model multiple substrates in future work.

Several aspects of the MEC system were overlooked in this study. We
did not consider rate hydrogen production directly since we do not
have a time series for either the hydrogen production rate or the
population of hydrogen consuming methanogens. If we were able to measure
either of those quantities, we could repeat this analysis to determine
precisely which parameters exert the most influence on the hydrogen
production rate itself. However, literature has shown that hydrogen
production directly correlates with current generation so the findings
presented in the study do represent the variation of hydrogen generation
from the MEC. The DAE model has the advantage of being computationally
inexpensive compared to PDE models, but the latter may provide a more
accurate characterization of the biofilm beyond merely concentration
in a well mixed compartment. Future work could analyze the bifurcations
in biofilm concentration in the context of a PDE model.

\section{Acknowledgements}

HJD and DMB are partially supported by NSF grant DMS-1225878. ZJR
and LL acknowledge support from NSF grant CBET-1510682 and ONR grant
N000141612210.  We thank Stephen Campbell for his help with questions about DAEs.

\bibliographystyle{siamplain}
\bibliography{MicrobialElectrolysisCell}

\begin{thebibliography}{10}

\bibitem{batstone_iwa_2002}
{\sc D.~J. Batstone, J.~Keller, I.~Angelidaki, S.~V. Kalyuzhnyi, S.~G.
  Pavlostathis, A.~Rozzi, W.~T.~M. Sanders, H.~Siegrist, and V.~A. Vavilin},
  {\em The {IWA} {Anaerobic} {Digestion} {Model} {No} 1 ({ADM}1)}, Water
  Science and Technology, 45 (2002), pp.~65--73.

\bibitem{beardmore_stability_1998}
{\sc R.~E. Beardmore}, {\em Stability and bifurcation properties of index-1
  {DAEs}}, Numerical Algorithms, 19 (1998), pp.~43--53.

\bibitem{brenan_numerical_1995}
{\sc K.~Brenan, S.~Campbell, and L.~Petzold}, {\em Numerical {Solution} of
  {Initial}-{Value} {Problems} in {Differential}-{Algebraic} {Equations}},
  Classics in {Applied} {Mathematics}, Society for Industrial and Applied
  Mathematics, Jan. 1995.

\bibitem{call_hydrogen_2008}
{\sc D.~Call and B.~E. Logan}, {\em Hydrogen production in a single chamber
  microbial electrolysis cell lacking a membrane}, Environmental Science \&
  Technology, 42 (2008), pp.~3401--3406.

\bibitem{chaudhuri_electricity_2003}
{\sc S.~K. Chaudhuri and D.~R. Lovley}, {\em Electricity generation by direct
  oxidation of glucose in mediatorless microbial fuel cells}, Nature
  Biotechnology, 21 (2003), pp.~1229--1232.

\bibitem{hansen_single-nutrient_1980}
{\sc S.~R. Hansen and S.~P. Hubbell}, {\em Single-nutrient microbial
  competition: qualitative agreement between experimental and theoretically
  forecast outcomes}, Science (New York, N.Y.), 207 (1980), pp.~1491--1493.

\bibitem{hindmarsh2005sundials}
{\sc A.~C. Hindmarsh, P.~N. Brown, K.~E. Grant, S.~L. Lee, R.~Serban, D.~E.
  Shumaker, and C.~S. Woodward}, {\em {SUNDIALS}: Suite of nonlinear and
  differential/algebraic equation solvers}, ACM Transactions on Mathematical
  Software (TOMS), 31 (2005), pp.~363--396.

\bibitem{hsu_mathematical_1977}
{\sc S.~Hsu, S.~Hubbell, and P.~Waltman}, {\em A {Mathematical} {Theory} for
  {Single}-{Nutrient} {Competition} in {Continuous} {Cultures} of
  {Micro}-{Organisms}}, SIAM Journal on Applied Mathematics, 32 (1977),
  pp.~366--383.

\bibitem{kato_marcus_conduction-based_2007}
{\sc A.~Kato~Marcus, C.~I. Torres, and B.~E. Rittmann}, {\em Conduction-based
  modeling of the biofilm anode of a microbial fuel cell}, Biotechnology and
  Bioengineering, 98 (2007), pp.~1171--1182.

\bibitem{liu_electrochemically_2005}
{\sc H.~Liu, S.~Grot, and B.~E. Logan}, {\em Electrochemically {Assisted}
  {Microbial} {Production} of {Hydrogen} from {Acetate}}, Environmental Science
  \& Technology, 39 (2005), pp.~4317--4320.

\bibitem{logan_exoelectrogenic_2009}
{\sc B.~E. Logan}, {\em Exoelectrogenic bacteria that power microbial fuel
  cells}, Nature Reviews Microbiology, 7 (2009), pp.~375--381.

\bibitem{logan_microbial_2006}
{\sc B.~E. Logan, B.~Hamelers, R.~Rozendal, U.~Schr\"{o}der, J.~Keller,
  S.~Freguia, P.~Aelterman, W.~Verstraete, and K.~Rabaey}, {\em Microbial fuel
  cells: Methodology and technology}, Environmental Science \& Technology, 40
  (2006), pp.~5181--5192.

\bibitem{lu_nickel_2016}
{\sc L.~Lu, D.~Hou, Y.~Fang, Y.~Huang, and Z.~J. Ren}, {\em Nickel based
  catalysts for highly efficient {H}2 evolution from wastewater in microbial
  electrolysis cells}, Electrochimica Acta, 206 (2016), pp.~381--387.

\bibitem{lu_active_2016}
{\sc L.~Lu, D.~Hou, X.~Wang, D.~Jassby, and Z.~J. Ren}, {\em Active {H}2
  {Harvesting} {Prevents} {Methanogenesis} in {Microbial} {Electrolysis}
  {Cells}}, Environmental Science \& Technology Letters, 3 (2016),
  pp.~286--290.

\bibitem{lu_hydrogen_2009}
{\sc L.~Lu, N.~Ren, D.~Xing, and B.~E. Logan}, {\em Hydrogen production with
  effluent from an ethanol-h$_2$-coproducing fermentation reactor using a
  single-chamber microbial electrolysis cell}, Biosensors and Bioelectronics,
  24 (2009), pp.~3055--3060.

\bibitem{lu_microbial_2016}
{\sc L.~Lu and Z.~J. Ren}, {\em Microbial electrolysis cells for waste
  biorefinery: {A} state of the art review}, Bioresource Technology, 215
  (2016), pp.~254--264.

\bibitem{meiss_differential_2017}
{\sc J.~Meiss}, {\em Differential {Dynamical} {Systems}, {Revised} {Edition}},
  Mathematical {Modeling} and {Computation}, Society for Industrial and Applied
  Mathematics, Jan. 2017.

\bibitem{monod_recherches_1942}
{\sc J.~Monod}, {\em Recherches sur la croissance des cultures
  bact\'{e}riennes}, Hermann, 1958.
\newblock Google-Books-ID: Gi9rAAAAMAAJ.

\bibitem{noren_clarifying_2005}
{\sc D.~A. Noren and M.~A. Hoffman}, {\em Clarifying the butler-volmer equation
  and related approximations for calculating activation losses in solid oxide
  fuel cell models}, Journal of Power Sources, 152 (2005), pp.~175--181.

\bibitem{novick_description_1950}
{\sc A.~Novick and L.~Szilard}, {\em Description of the {Chemostat}}, Science,
  112 (1950), pp.~715--716.

\bibitem{oliveira_overview_2013}
{\sc V.~B. Oliveira, M.~Sim\~{o}es, L.~F. Melo, and A.~M. F.~R. Pinto}, {\em
  Overview on the developments of microbial fuel cells}, Biochemical
  Engineering Journal, 73 (2013), pp.~53--64.

\bibitem{ortiz-martinez_developments_2015}
{\sc V.~M. Ortiz-Mart\'{i}­nez, M.~J. Salar-Garc\'{i}­a, A.~P. de~los
  R\'{i}­­os, F.~J. Hern\'{a}­ndez-Fern\'{a}­ndez, J.~A. Egea, and L.~J.
  Lozano}, {\em Developments in microbial fuel cell modeling}, Chemical
  Engineering Journal, 271 (2015), pp.~50--60.

\bibitem{perko_differential_2001}
{\sc L.~Perko}, {\em Differential {Equations} and {Dynamical} {Systems}},
  Springer, 2001.

\bibitem{picioreanu_computational_2007}
{\sc C.~Picioreanu, I.~M. Head, K.~P. Katuri, M.~C.~M. van Loosdrecht, and
  K.~Scott}, {\em A computational model for biofilm-based microbial fuel
  cells}, Water Research, 41 (2007), pp.~2921--2940.

\bibitem{picioreanu_model_2010}
{\sc C.~Picioreanu, M.~C.~M. van Loosdrecht, T.~P. Curtis, and K.~Scott}, {\em
  Model based evaluation of the effect of {pH} and electrode geometry on
  microbial fuel cell performance}, Bioelectrochemistry, 78 (2010), pp.~8--24.

\bibitem{picioreanu_mathematical_2008}
{\sc C.~Picioreanu, M.~C.~M. van Loosdrecht, K.~P. Katuri, K.~Scott, and I.~M.
  Head}, {\em Mathematical model for microbial fuel cells with anodic biofilms
  and anaerobic digestion}, Water Science and Technology: A Journal of the
  International Association on Water Pollution Research, 57 (2008),
  pp.~965--971.

\bibitem{pinto_two-population_2010}
{\sc R.~Pinto, B.~Srinivasan, M.-F. Manuel, and B.~Tartakovsky}, {\em A
  two-population bio-electrochemical model of a microbial fuel cell},
  Bioresource Technology, 101 (2010), pp.~5256--5265.

\bibitem{pinto2011}
{\sc R.~P. Pinto, B.~Srinivasan, A.~Escapa, and B.~Tartakovsky}, {\em
  Multi-population model of a microbial electrolysis cell}, Environmental
  Science \& Technology, 45 (2011), pp.~5039--5046.

\bibitem{rabaey_microbial_2010}
{\sc K.~Rabaey and R.~Rozendal}, {\em Microbial electrosynthesis - revisiting
  the electrical route for microbial production}, Nature Reviews Microbiology,
  8 (2010), pp.~706--716.

\bibitem{ciarlet_1990}
{\sc P.~J. Rabier and W.~C. Rheinboldt}, {\em Theoretical and numerical
  analysis of differential-algebraic equations}, in Handbook of {Numerical}
  {Analysis}, P.~Ciarlet and J.~Lions, eds., no.~v. 8 in Handbook of
  {Numerical} {Analysis}, North-Holland, 1990.

\bibitem{recio-garrido_modeling_2016}
{\sc D.~Recio-Garrido, M.~Perrier, and B.~Tartakovsky}, {\em Modeling,
  optimization and control of bioelectrochemical systems}, Chemical Engineering
  Journal, 289 (2016), pp.~180--190.

\bibitem{reich_local_1995}
{\sc S.~Reich}, {\em On the local qualitative behavior of
  differential-algebraic equations}, Circuits, Systems and Signal Processing,
  14 (1995), pp.~427--443.

\bibitem{riaza_stability_2002}
{\sc R.~Riaza}, {\em Stability {Issues} in {Regular} and {Noncritical}
  {Singular} {DAEs}}, Acta Applicandae Mathematica, 73 (2002), pp.~301--336.

\bibitem{riaza_2008}
{\sc R.~Riaza}, {\em Differential-{Algebraic} {Systems}: {Analytical} {Aspects}
  and {Circuit} {Applications}}, World Scientific, 2008.

\bibitem{rozendal_principle_2006}
{\sc R.~A. Rozendal, H.~V.~M. Hamelers, G.~J.~W. Euverink, S.~J. Metz, and
  C.~J.~N. Buisman}, {\em Principle and perspectives of hydrogen production
  through biocatalyzed electrolysis}, International Journal of Hydrogen Energy,
  31 (2006), pp.~1632--1640.

\bibitem{smith_theory_1995}
{\sc H.~L. Smith and P.~Waltman}, {\em The {Theory} of the {Chemostat}:
  {Dynamics} of {Microbial} {Competition}}, Cambridge University Press, Jan.
  1995.
\newblock Google-Books-ID: wFLdVo89vq8C.

\bibitem{sotomayor_generic_1973}
{\sc J.~Sotomayor}, {\em Generic {Bifurcations} of {Dynamical} {Systems}}, in
  Dynamical {Systems}, M.~M. Peixoto, ed., Academic Press, 1973, pp.~561--582.

\bibitem{venkatasubramanian_local_1995}
{\sc V.~Venkatasubramanian, H.~Schattler, and J.~Zaborszky}, {\em Local
  bifurcations and feasibility regions in differential-algebraic systems}, IEEE
  Transactions on Automatic Control, 40 (1995), pp.~1992--2013.

\bibitem{wang_comprehensive_2013}
{\sc H.~Wang and Z.~J. Ren}, {\em A comprehensive review of microbial
  electrochemical systems as a platform technology}, Biotechnology Advances, 31
  (2013), pp.~1796--1807.

\bibitem{zhang_modelling_1995}
{\sc X.-C. Zhang and A.~Halme}, {\em Modelling of a microbial fuel cell
  process}, Biotechnology Letters, 17 (1995), pp.~809--814.

\end{thebibliography}

\end{document}